\documentclass[11pt]{article}

\usepackage{amssymb}
\usepackage{amsmath}
\usepackage{amsthm}
\usepackage{enumerate}

\newtheoremstyle{plain}%
    {8pt plus2pt minus4pt}%
    {8pt plus2pt minus4pt}%
    {\itshape}%
    {}%
    {\bfseries\scshape}%
    {}%
    {6pt}
    {}%

\newtheoremstyle{remark}%
    {8pt plus2pt minus4pt}%
    {8pt plus2pt minus4pt}%
    {\upshape}
    {}%
    {\bfseries\scshape}%
    {}%
    {6pt}
    {}%

\theoremstyle{plain}
\newtheorem{thm}{Theorem}
\newtheorem{cor}[thm]{Corollary}

\newtheorem{conj}[thm]{Conjecture}

\newtheorem{defn}[thm]{Definition}

\theoremstyle{remark}

\newcommand{\vv}{\vec{v}}
\newcommand{\ww}{\vec{w}}
\newcommand{\xx}{\vec{x}}
\newcommand{\yy}{\vec{y}}
\newcommand{\uu}{\vec{u}}

\title{Coloring the cube with rainbow cycles}
\author{Dhruv
Mubayi\thanks{Department of Mathematics, Statistics, and Computer
Science, University of Illinois, Chicago, IL 60607;  email:
mubayi@math.uic.edu; research  supported in part by  NSF grant DMS
0969092} \quad \quad Randall Stading\thanks{Department of Mathematics, Statistics, and Computer
Science, University of Illinois, Chicago, IL 60607;  email: rstadi1@uic.edu}}

\usepackage{tikz}
\usepackage{amsfonts}
\usepackage{amssymb}
\begin{document}

\maketitle

\begin{abstract}
For every even positive integer $k\ge 4$ let $f(n,k)$ denote the minimim number of colors required to color the edges of the $n$-dimensional cube $Q_n$, so that the edges of every copy of $k$-cycle $C_k$  receive $k$ distinct colors. 
Faudree, Gy\'arf\'as, Lesniak and Schelp proved that $f(n,4)=n$ for $n=4$ or $n>5$. We consider larger $k$ and prove that if $k \equiv 0$ (mod 4), then there are positive constants $c_1, c_2$ depending only on $k$ such that
$$c_1n^{k/4} < f(n,k) < c_2 n^{k/4}.$$
Our upper bound uses an old construction of Bose and Chowla of generalized Sidon sets.
For  $k \equiv 2$ (mod 4), the situation seems more complicated. For  the smallest case $k=6$ we show that $$n \le f(n, 6) < n^{1+o(1)}.$$ The upper bound is obtained from  Behrend's construction of a subset of the integers with no three term arithmetic progression.
\end{abstract}
\vspace{25pt}
\section{Introduction}
Given graphs $G$ and $H$, and an integer $q\le |E(H)|$, a $(G,H,q)$-coloring is an
edge-coloring of $G$ such that the edges of every copy of $H$ in $G$ receive at least $q$ colors. Let
$f(G, H, q)$ be the minimum number of colors in a $(G,H,q)$-coloring. This general problem is hopeless in most cases, for example, when $G$ and $H$ are cliques, and $q=2$, determining it is equivalent to determining the multicolor Ramsey number $R_k(p)$ which is a longstanding open problem. There has been more success in determining $f(G,H,q)$  when $G$ and $H$ are not cliques or when $q>2$ (or both). Many Ramsey problems have received considerable attention when studied on the $n$-dimensional cube. The papers
\cite{ARSV, AHMS} are examples where anti-Ramsey problems for subcubes in cubes and problems about monochromatic cycles in cubes are investigated. In \cite{DO}, Offner found the exact value for the maximum number of colors for which it is possible to edge color the hypercube so that all subcubes of dimension $d$ contain all colors. Related Tur\'an type problems for subcubes in cubes have been studied in \cite{AKS}.

Rainbow cycles have also been well studied as subgraphs of $K_n$. Erd\H os, Simonovits and S\'os \cite{ESS} introduced $AR(n,H)$, the maximum number of colors in an edge coloring of $K_n$ such that it contains no rainbow copy of $H$, and provided a conjecture when $H$ is a cycle and showed that their conjecture was true when  $H = C_3$. Alon \cite{NA} proved their conjecture for cycles of length four  and Montellano-Ballesteros and Neumann-Lara \cite{MBNL} proved the conjecture for all cycles in 2003. More recently, Choi \cite{JC} gave a shorter proof of the conjecture.

We continue this theme in the current note and let $G=Q_n$, the $n$-dimensional cube,  and $H=C_k$, the cycle of length $k$.  Our focus is on  $q=|E(H)|$, in which case we will call a $(G,H,q)$-coloring an $H$-rainbow coloring, assuming that $G$ is obvious from context (in this paper $G=Q_n$ always).

\begin{defn} For $4\le k \le 2^n$, let $f(n,k)=f(Q_n, C_k, k)$
be the minimum number of colors in a $C_k$-rainbow coloring of $Q_n$. \end{defn}

The smallest case $f(n,4)$  was studied by Faudree, Gy\'arf\'as, Lesniak and Schelp~\cite{FGLS} who proved that the trivial lower bound of $n$ is tight by providing, for all $n \ge 6$,  a $C_4$-rainbow coloring with $n$ colors. We consider larger $k$. Our first result determines the order of magnitude of $f(n,k)$  for $k \equiv 0$ (mod 4).

\begin{thm} \label{mod4}  Fix a positive $k \equiv 0$ (mod 4).  There are constants $c_1, c_2>0$ depending only on $k$ such that
$$c_1n^{k/4} < f(n,k) < c_2 n^{k/4}.$$
\end{thm}

The case $k\equiv 2$ (mod 4) seems more complicated. Our results imply that for such fixed $k$
there are positive constants $c_1', c_2'$ with 
$$c_1' n^{\lfloor k/4 \rfloor} < f(n,k) < c_2' n^{\lceil k/4 \rceil}.$$
We believe that the lower bound is closer to the truth.  As evidence for this, we tackle the smallest case in this range, $k=6$.  As we will observe later, the lower bound $f(n,6) \ge n$ is trivial for $n \ge 3$, and we obtain the following upper bound.

\begin{thm} \label{6} For every $\epsilon > 0$ there exists $n_0$ such that for $n>n_0$ we have~$f(n,6) \leq n^{1 + \epsilon}$.
\end{thm}

It is rather easy to see that $f(Q_n,Q_3,12) = f(Q_n,C_6,3)$.
Indeed, if $Q_n$ is edge-colored so that every $Q_3$ is rainbow, then every $C_6$ is rainbow since each one is contained in a rainbow $Q_3$ and so $f(Q_n,Q_3,12) \geq f(Q_n,C_6,3)$. On the other hand, it is easy to see that any two edges of a $Q_3$ lie in some $C_6$ and therefore if $Q_n$ is edge-colored so that every  $C_6$ is rainbow then every $Q_3$ must also be rainbow and so $f(Q_n,Q_3,12) \leq f(Q_n,C_6,3)$. 

Since $C_4=Q_2$, the following corollary can also be considered an extension of the result \cite{FGLS} to subcubes.

\begin{cor} As $n \rightarrow \infty$, we have $f(Q_n,Q_3,12) =n^{1+o(1)}$.
\end{cor}

We will consider the vertices of  $Q_n$  as binary vectors of length $n$ or as subsets of $[n]=\{1, \ldots, n\}$, depending on the context (with the natural bijection $\vv \leftrightarrow v$ where $\vv$ is the incidence vector for $v \subset[n]$, i.e. $\vv_i=1$ iff $i \in v$). 
 In particular, whenever we write $v-w$ we mean set theoretic difference, $v \cup w$ or $v \cap w$ we mean set union/intersection and when we write $\vv \pm \ww$ we mean 
 vector addition/subtraction modulo 2.
  We  write $e_i$ for the standard basis vector, so $e_i$ is one in the $i$th coordinate and zero in all other coordinates.  Given an edge $f=uv$ of $Q_n$ where $\vv=\uu +e_s$ for some $s$, we say that $v$ is the top vertex of $f$ and $u$ is the bottom vertex.
We will say the an edge is on level $i$ of $Q_n$ if its bottom vertex corresponds to a vector with $i-1$ ones and the top vertex to a vector with $i$ ones.

\section{Proof of Theorem \ref{mod4}}
The lower bound in Theorem \ref{mod4} follows from the easy observation that in a $C_k$-rainbow coloring all edges at level $k/4$ must receive distinct colors. Indeed, given any two such edges $f_1=vw$ and $f_2=xy$, where $\ww=\vv+e_i$ and $\yy=\xx+e_j$, it suffices to find a copy of $C_k$ containing $f_1$ and $f_2$. If $f_1$ and $f_2$ are incident then it is clear that we can find a $C_k$ containing them as long as $n>k$ which we may clearly assume. The two cases are illustrated below where $r = k/2 - 2$ and $s_i \notin w \cup y$ for all $i \in \lbrace 1,...,r \rbrace$.
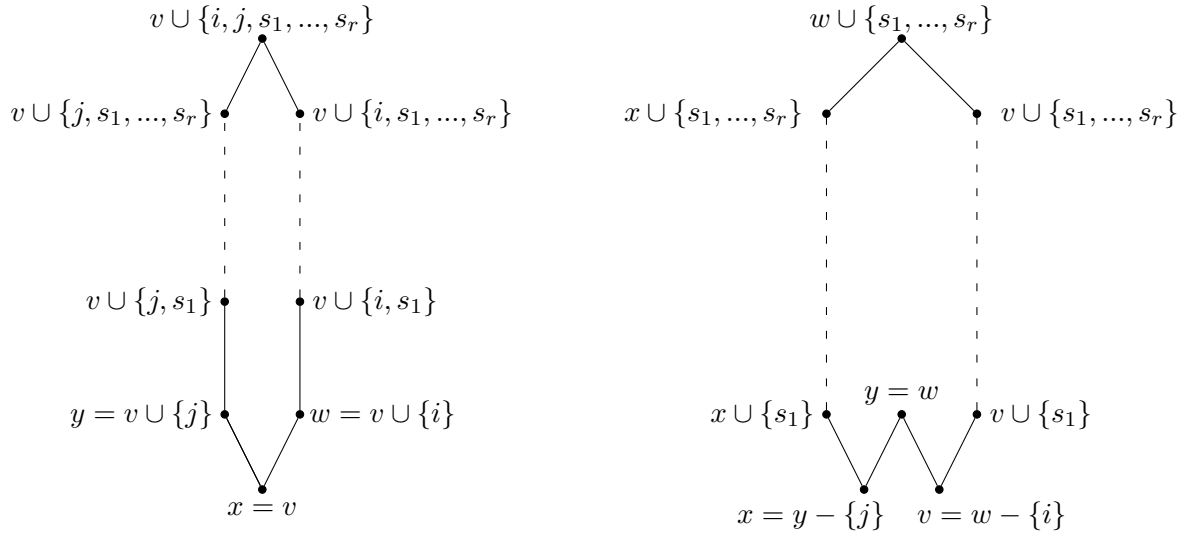
  \begin{center}
  \begin{tikzpicture}[auto]
  \tikzstyle{vertex}=[circle,draw = black,fill=black,minimum size=1pt,inner sep=1pt]
  \node[vertex] (V) at (-3.5,0) {};
  \draw (-3.5,-.25) node {$x=v$};
  \node[vertex] (VL1) at (-4,1) {};
  \draw (-5.1,1)node {$y=v\cup\lbrace j \rbrace$};
  \node[vertex] (VR1) at (-3,1) {};
  \draw (-1.9,1) node {$w=v\cup\lbrace i \rbrace$};
  \node[vertex] (VL2) at (-4,2.5) {};
  \draw (-5,2.5) node {$v\cup\lbrace j,s_1 \rbrace$};
  \node[vertex] (VR2) at (-3,2.5) {};
  \draw (-2,2.5) node {$v\cup\lbrace i,s_1 \rbrace$};
  \node[vertex] (VL3) at (-4,5) {};
  \draw (-5.5,5) node {$v\cup\lbrace j,s_1,...,s_r \rbrace$};
  \node[vertex] (VR3) at (-3,5) {};
  \draw (-1.5,5) node {$v\cup\lbrace i,s_1,...,s_r \rbrace$};
  \node[vertex] (VT) at (-3.5,6) {};
  \draw (-3.5,6.25) node {$v\cup\lbrace i,j,s_1,...,s_r \rbrace$};
  \draw (VL1) to node {} (V);
  \draw (VR1) to node {} (V);
  \draw (VR1) to node {} (VR2);
  \draw (VL1) to node {} (V);
  \draw (VL1) to node {} (VL2);
  \draw [loosely dashed](VR2) to node {} (VR3);
  \draw [loosely dashed](VL2) to node {} (VL3);
  \draw (VR3) to node {} (VT);
  \draw (VL3) to node {} (VT);
  \node[vertex] (W) at (5,1) {};
  \draw (5,1.25) node {$y=w$};
  \node[vertex] (WL1) at (4.5,0) {};
  \draw (3.8,-.35) node {$x=y-\lbrace j \rbrace$};
  \node[vertex] (WR1) at (5.5,0) {};
  \draw (6.2,-.35) node {$v=w-\lbrace i \rbrace$};
  \node[vertex] (WL2) at (4,1) {};
  \draw (3.15,1) node {$x\cup\lbrace s_1 \rbrace$};
  \node[vertex] (WR2) at (6,1) {};
  \draw (6.85,1) node {$v\cup\lbrace s_1 \rbrace$};
  \node[vertex] (WL3) at (4,5) {};
  \draw (2.5,5) node {$x\cup\lbrace s_1,...,s_{r} \rbrace$};
  \node[vertex] (WR3) at (6,5) {};
  \draw (7.5,5) node {$v\cup\lbrace s_1,...,s_{r} \rbrace$};
  \node[vertex] (WT) at (5,6) {};
  \draw (5,6.25) node {$w\cup\lbrace s_1,...,s_{r} \rbrace$};
  \draw (W) to node {} (WL1);
  \draw (WL1) to node {} (WL2);
  \draw [loosely dashed](WL2) to node {} (WL3);
  \draw (WL3) to node {} (WT);
  \draw (W) to node {} (WR1);
  \draw (WR1) to node {} (WR2);
  \draw [loosely dashed](WR2) to node {} (WR3);
  \draw (WR3) to node {} (WT);
  \end{tikzpicture}
  \end{center}
Now, suppose $f_1$ and $f_2$ are not incident. We know that $|x \triangle v| \leq k/2 - 2$ since $x$ and $v$ are each sets of size $k/4 - 1$. By successively deleting elements of $v$ and $x$ in the appropriate order, we can obtain a $v,x$-path of length $k/2 - 2$. Then, since $w$ and $y$ are sets of size $k/4$, we may find a $w,y$-path of length $k/2$ between them by successively adding the elements of $y$ to $w$ and vice versa along with extra elements as needed. The two paths along with the edges $vw$ and $xy$ form a cycle of length $k$. This is shown in the following diagram. Let $y-w = \lbrace y_1,...,y_m \rbrace$, $w-y = \lbrace w_1,...,w_m\rbrace$ and $w\cap y = \lbrace z_1,...,z_l\rbrace$ where $m + l = k/4$. Let $\lbrace s_1,...,s_r \rbrace$ again be a set such that $s_i \notin y \cup w$ with $r = k/4 - m$.
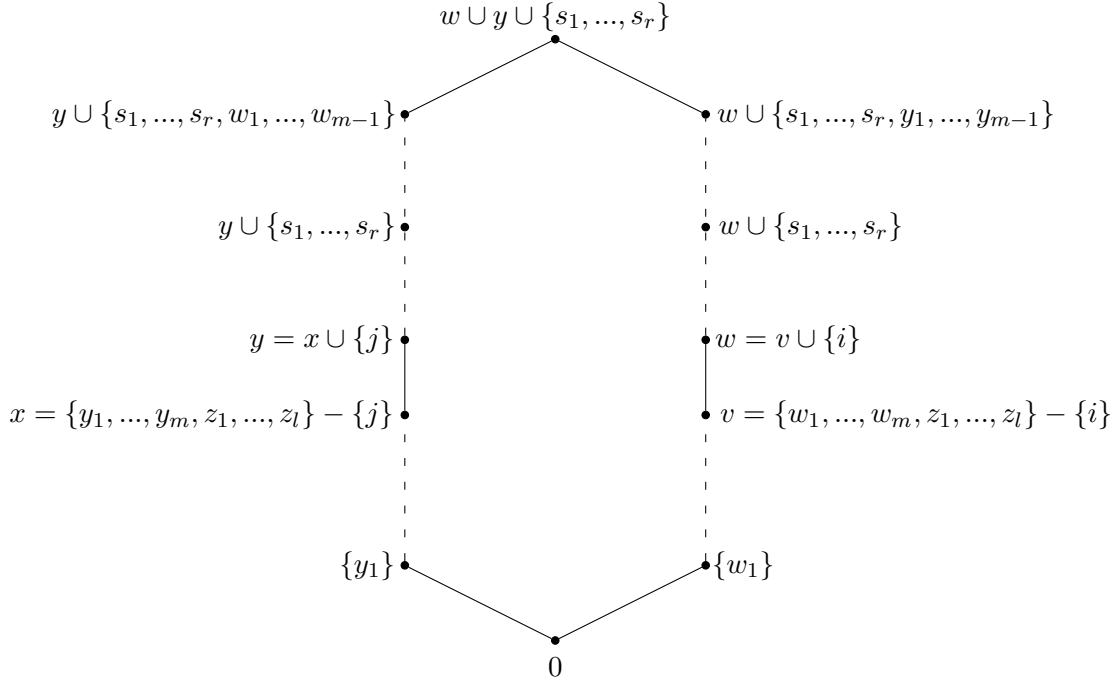
  \begin{center}
  \begin{tikzpicture}[auto]
  \tikzstyle{vertex}=[circle,draw = black,fill=black,minimum size=1pt,inner sep=1pt]
  \node [vertex] (G) at (0,2) {};
  \draw (0,1.65) node {$0$};
  \node [vertex] (GL1) at (-2,3) {};
  \draw (-2.5,3) node {$\lbrace y_1 \rbrace$};
  \node [vertex] (GR1) at (2,3) {};
  \draw (2.5,3) node {$\lbrace w_1 \rbrace$};
  \node [vertex] (GL3) at (-2,5) {};
  \draw (-4.7,5) node {$x=\lbrace y_1,...,y_m,z_1,...,z_l \rbrace - \lbrace j \rbrace$};
  \node [vertex] (GR3) at (2,5) {};
  \draw (4.8,5) node {$v=\lbrace w_1,...,w_m,z_1,...,z_l \rbrace - \lbrace i \rbrace$};
  \node [vertex] (GL4) at (-2,6) {};
  \draw (-3.1,6) node {$y=x \cup \lbrace j \rbrace$};
  \node [vertex] (GR4) at (2,6) {};
  \draw (3.1,6) node {$w=v \cup \lbrace i \rbrace$};
  \node [vertex] (GL5) at (-2,7.5) {};
  \draw (-3.3,7.5) node {$y \cup \lbrace {s_1,...,s_r} \rbrace$};
  \node [vertex] (GR5) at (2,7.5) {};
  \draw (3.4,7.5) node {$w \cup \lbrace {s_1,...,s_r} \rbrace$};
  \node [vertex] (GL6) at (-2,9) {};
  \draw (-4.4,9) node {$y \cup \lbrace {s_1,...,s_r,w_1,...,w_{m-1}} \rbrace$};
  \node [vertex] (GR6) at (2,9) {};
  \draw (4.4,9) node {$w \cup \lbrace {s_1,...,s_r,y_1,...,y_{m-1}} \rbrace$};
  \node [vertex] (GT) at (0,10) {};
  \draw (0,10.3) node {$w \cup y \cup \lbrace {s_1,...,s_r} \rbrace$};
  \draw (GL1) to (G) to (GR1);
  \draw (GL6) to (GT) to (GR6);
  \draw (GL3) to (GL4);
  \draw (GR3) to (GR4); 
  \draw [loosely dashed] (GL1) to (GL3);
  \draw [loosely dashed] (GR1) to (GR3);
  \draw [loosely dashed] (GL4) to (GL6);
  \draw [loosely dashed] (GR4) to (GR6);
  \end{tikzpicture}
  \end{center}
For the upper bound we need a classical construction of generalized Sidon sets by Bose and Chowla.  A $B_t$-set $S=\{s_1, \ldots, s_n\}$ is a set of integers such that if $1\le i_1 \le i_2 \le \cdots \le i_{t}\le n$ and $1\le j_1 \le j_2 \le \cdots \le j_{t}\le n$, then
$$ s_{i_1}+\cdots +s_{i_t} \neq s_{j_1}+\cdots +s_{j_t}
$$
unless $(i_1, \ldots, i_t)=(j_1, \ldots, j_t)$.  A consequence of this is that if $P,Q$ are nonempty disjoint subsets of $[n]$ with $|P|=|Q|\le t$, then
\begin{equation} \label{sidon}
\sum_{i \in P} s_i \neq \sum_{j \in Q} s_j.
\end{equation}

 The result below is phrased in a form that is suitable for our use later. 

\begin{thm} {\bf (Bose-Chowla \cite{BC})} \label{bc}
For each fixed $t \ge 2$, there is a constant $A>1$ such that for all $n$, there is a $B_t$-set $S=\{s_1, \ldots, s_n\} \subset \{1, 2, \ldots, \lfloor An^t \rfloor\}$.
\end{thm}

Now we provide the upper bound construction for Theorem \ref{mod4}.

{\bf Construction 1.} Let $t=k/4-1$ and $S=\{s_1, \ldots, s_n\} \subset \{1, 2, \ldots, \lfloor An^{t} \rfloor\}$ be a $B_{t}$-set as above.   For each $v \in V(Q_n)$, let 
$$a(v)= \sum_{i=1}^n \vv_i s_i= \sum_{i: \vv_i=1} s_i.$$ 
Given $vw \in E(Q_n)$
with $\ww=\vv+e_j$,   let $M=\lceil kAn^{t}\rceil$, and let
$$d(vw)= a(v) + Mj.$$
Suppose further that $vw$ is at level $p$ and $p'$ is the congruence class of $p$ modulo $k/2$.  Then the color of the edge $vw$ is 
$$\chi(vw)= (d(vw), p').  \qed $$

Let us now argue that this construction yields the upper bound in Theorem \ref{mod4}.  First, the number of colors is at most 
$$\max_{vw} d(vw) \times \frac{k}{2}\le (n \cdot \max s_i + Mn)\frac{k}{2}\le \frac{nk}{2}An^t + \frac{nk}{2}M < k^2An^{t+1}= k^2An^{k/4}$$ as desired.  Now we show that this is a $C_k$-rainbow coloring. Suppose for contradiction that $H$ is a copy of $C_k$ in $Q_n$ and $f_1=vw, f_2=xy$ are distinct edges of $H$ with $\chi(f_1)=\chi(f_2)$. Since $H$ spans at most $k/2$ levels,  $f_1$ and $f_2$ cannot lie in levels that differ by more than $k/2$, so $\chi(f_1)\neq 
\chi(f_2)$ unless $f_1$ and $f_2$ are in the same level which we may henceforth assume.  Let $v,x$ be the bottom vertices of $f_1, f_2$, and $\ww=\vv+e_i, \yy=\xx+e_j$.  Assume without loss of generality that $i \le j$. If $v=x$, then $$a(v)+Mi=d(vw)=d(xy)=a(x)+Mj=a(v)+Mj.$$ This  implies that $i=j$ and
 contradicts the fact that $f_1 \neq f_2$.  We may therefore assume that $v\neq x$. 
Similarly, if $w=y$, then $i < j$ and $$a(w)-s_i+Mi=a(v)+Mi = d(vw)=d(xy)= a(x)+Mj= a(y)-s_j+Mj=a(w)-s_j+Mj.$$ This implies the contradiction   $s_j-s_i=M(j-i)\ge M > An^t > s_j-s_i$.
Consequently, we may assume that $vw$ and $xy$ share no vertex. 
 
  If $|v \triangle x| >k/2$, then any $v,x$-path in $Q_n$ has length more than $k/2$ so there can be no cycle of length $k$ containing both $v$ and $x$, contradiction. So we may assume that $|v \triangle x|\le k/2$. Now $\chi(vw)=\chi(xy)$ implies that $$a(v)+Mi=d(vw)=d(xy)=a(x)+Mj$$ and this yields
 $$M(j-i)=Mj-Mi= a(v)-a(x)= a(v-x)-a(x-v) \le \frac{|v \triangle x|}{2} An^t\le \frac{k}{4}An^{t}<M.$$
 Consequently, we may assume that $i=j$, $a(v)=a(x)$,
 $a(v-x)=a(x-v)$ and $|v \triangle x|=|w \triangle y|$.
  If $|v-x|=|x-v|\le k/4-1$, then
 $$a(v-x) = \sum_{i\in v-x} s_i \neq  \sum_{j \in x-v} s_j= a(x-v)$$ due to (\ref{sidon}), the definition of $S$ and $t=k/4-1$.  So we may assume that $|v-x|=|x-v|= k/4$ and 
 $|w \triangle y|=|v \triangle x|=k/2$. This implies that dist$_{Q_n}(w,y)=$ dist$_{Q_n}(v,x)=k/2$. Together with edges $f_1, f_2$, we conclude that $C$ must have at least $k+2$ edges, contradiction. \qed

\section{Proof of Theorem \ref{6}}

We will first show the lower bound $f(n,6) \ge n$ for $n \ge 3$. It is immediate that no two edges incident with 0 receive the same color, for otherwise there would be two edges of the same color on level 1 of $Q_n$ which we could easily extend to a non-rainbow $C_6$.  Indeed,  let $i,j,k$ be distinct and consider the following $C_6$: 
$$0 \quad e_i \quad e_i+e_k \quad e_k\quad  e_j+e_k\quad e_j \quad 0.$$

To obtain the upper bound, we will give an explicit coloring that makes use of a classical construction of Behrend on sets of integers with no  arithmetic progression of size three. Let $r_3(N)$ denote the maximum size of a subset of $\{1,\ldots, N\}$ that contains no 3-term arithmetic progression.

\begin{thm} {\bf (Behrend \cite{B})}  \label{b}
 There is a $c>0$ such that if $N$ is sufficiently large, then
 $$r_3(N) > N^{1 - {c \over {\sqrt{\log N}}}}.$$
\end{thm} 
 
Theorem \ref{b} clearly implies that for $\epsilon > 0$ and sufficiently large $N$ we have $r_3(N) > N^{1 - \epsilon}$.  The error term $\epsilon$  was improved  recently by Elkin~\cite{E} (see \cite{GW} for a simpler proof) and using Elkin's result would give corresponding improvements in our result. 

{\bf Construction 2.}
Let $\epsilon > 0$ and $n$ be sufficiently large. Put $N=\lceil n^{1+\epsilon} \rceil$ and let  $S = \{s_1,...,s_n\}\subset \{1, \ldots, N\}$  contain no 3-term arithmetic progression.   Such a set exists by Theorem \ref{b} since 
 $$n > n^{1-\epsilon^2} = n^{(1-\epsilon)(1+\epsilon)}
 > N^{1-2\epsilon}.$$
  Let 
 $$a(v) = {\sum\limits_{i = 1}^{n} \vv_i\,s_i}.$$
 Consider the edge $vw$, where $\ww=\vv+e_k$. Let
 $$d(vw)= a(v) +2s_k \in Z_{2N}.$$
 We emphasize here that we are computing $d(vw)$ modulo $2N$.
 Suppose further that $vw$ is at level $p$ and $p'$ is the congruence class of $p$ modulo $3$.  Then the color of the edge $vw$ is 
$$\chi(vw)= (d(vw), p').  \qed $$
  
The number of colors used is at most $6N<n^{1+2\epsilon}$ as required.  We will now show that this is a $C_6$-rainbow coloring. Due to the second coordinate, it suffices to show that any
 two edges $f_1, f_2$ of a $C_6$ which are on the same of level of $Q_n$ receive different colors. If $f_1$ and $f_2$ are incident, then they meet either at their top vertices or bottom vertices. 
  
  If incident at their bottom vertices, the edges are colored as follows and thus are distinctly colored:
  
  
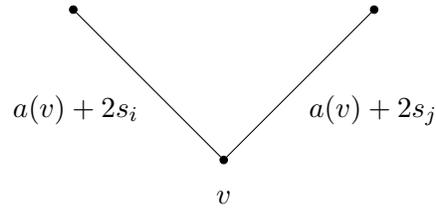
\begin{center}
  \begin{tikzpicture}[auto]
  \tikzstyle{vertex}=[circle,draw = black,fill=black,minimum size=1pt,inner sep=1pt]
  \node[vertex] (V) at (0,0) {};
  \node[vertex] (VL) at (-2,2) {};
  \node[vertex] (VR) at (2,2) {};
  \draw(0,-.5) node {$v$};
  \draw (V) to node {$a(v) + 2s_i$} (VL);
  \draw (V) to node [swap] {$a(v) + 2s_j$} (VR);
  \end{tikzpicture}
  \end{center}
  
  If incident at their top vertices, the edges lie on a $C_4$ and are therefore distinctly colored.
  
  
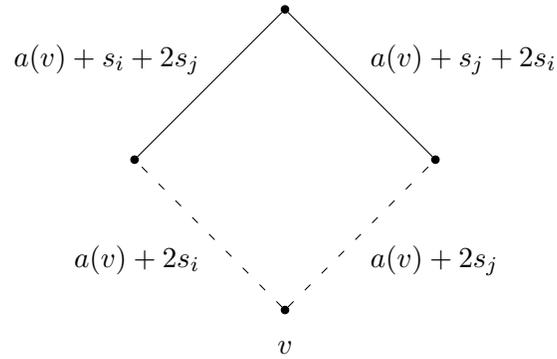
\begin{center}
  \begin{tikzpicture}[auto]
  \tikzstyle{vertex}=[circle,draw = black,fill=black,minimum size=1pt,inner sep=1pt]
  \node[vertex] (V) at (0,0) {};
  \node[vertex] (VL) at (-2,2) {};
  \node[vertex] (VR) at (2,2) {};
  \node[vertex] (VT) at (0,4) {};
  \draw(0,-.5) node {$v$};
  \draw [loosely dashed](V) to node {$a(v) + 2s_i$} (VL);
  \draw [loosely dashed](V) to node [swap] {$a(v) + 2s_j$} (VR);
  \draw (VR) to node [swap]{$a(v) + s_j + 2s_i$} (VT);
  \draw (VL) to node {$a(v) + s_i + 2s_j$} (VT);
  \end{tikzpicture}
  \end{center}
  
  If $f_1$ and $f_2$ are not incident, then there must be a path of length two between their bottom vertices. For if not, then they could not lie on a $C_6$ as the shortest path between their top vertices has length at least two. Moreover, the top vertices of $f_1$ and $f_2$ have symmetric difference precisely two since there is a path of length two between them.
  With these conditions, there are three ways the edges may be colored.
  
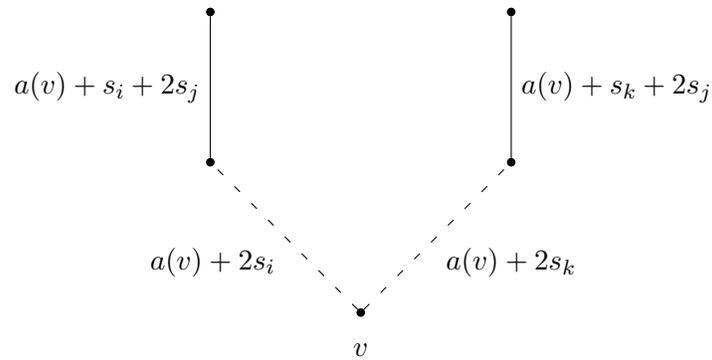
\begin{center}
  \begin{tikzpicture}[auto]
  \tikzstyle{vertex}=[circle,draw = black,fill=black,minimum size=1pt,inner sep=1pt]
  \node[vertex] (V) at (0,0) {};
  \node[vertex] (L1) at (-2,2) {};
  \node[vertex] (L2) at (-2,4) {};
  \node[vertex] (R1) at (2,2) {};
  \node[vertex] (R2) at (2,4) {};
  \draw [loosely dashed] (V) to node {$a(v) + 2s_i$} (L1);
  \draw (L1) to node {$a(v)+ s_i + 2s_j$}(L2);
  \draw [loosely dashed] (V) to node [swap]{$a(v) + 2s_k$} (R1);
  \draw (R1) to node [swap]{$a(v)+ s_k + 2s_j$} (R2);
  \draw(0,-.5) node {$v$};
  \end{tikzpicture}
  \end{center}
  
  
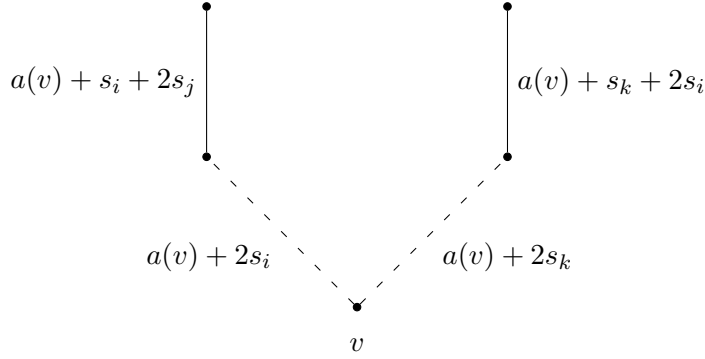
\begin{center}
  \begin{tikzpicture}[auto]
  \tikzstyle{vertex}=[circle,draw = black,fill=black,minimum size=1pt,inner sep=1pt]
  \node[vertex] (V) at (0,0) {};
  \node[vertex] (L1) at (-2,2) {};
  \node[vertex] (L2) at (-2,4) {};
  \node[vertex] (R1) at (2,2) {};
  \node[vertex] (R2) at (2,4) {};
  \draw [loosely dashed] (V) to node {$a(v) + 2s_i$} (L1);
  \draw (L1) to node {$a(v)+ s_i + 2s_j$}(L2);
  \draw [loosely dashed] (V) to node [swap]{$a(v) + 2s_k$} (R1);
  \draw (R1) to node [swap]{$a(v)+ s_k + 2s_i$} (R2);
  \draw(0,-.5) node {$v$};
  \end{tikzpicture}
  \end{center}
  
  
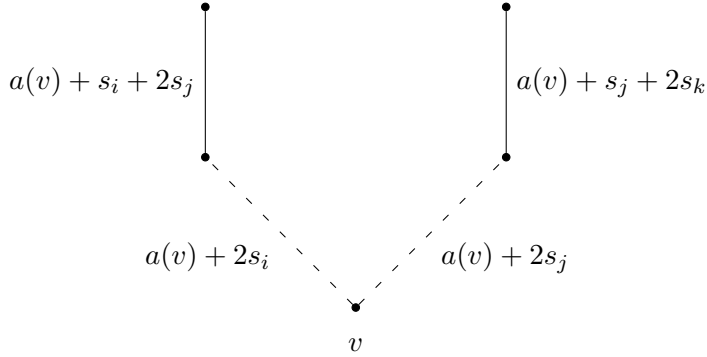
\begin{center}
  \begin{tikzpicture}[auto]
  \tikzstyle{vertex}=[circle,draw = black,fill=black,minimum size=1pt,inner sep=1pt]
  \node[vertex] (V) at (0,0) {};
  \node[vertex] (L1) at (-2,2) {};
  \node[vertex] (L2) at (-2,4) {};
  \node[vertex] (R1) at (2,2) {};
  \node[vertex] (R2) at (2,4) {};
  \draw [loosely dashed] (V) to node {$a(v) + 2s_i$} (L1);
  \draw (L1) to node {$a(v)+ s_i + 2s_j$}(L2);
  \draw [loosely dashed] (V) to node [swap]{$a(v) + 2s_j$} (R1);
  \draw (R1) to node [swap]{$a(v)+ s_j + 2s_k$} (R2);
  \draw(0,-.5) node {$v$};
  \end{tikzpicture}
  \end{center}
  
  In the first coloring, $s_i + 2s_j \neq s_k + 2s_j$ holds due to $i$ and $k$ being distinct. In the second and third colorings, $s_i + 2s_j \neq s_k + 2s_i$ and $s_i + 2s_j \neq s_j + 2s_k$ hold due to our set $S$ being free of three term arithmetic progressions.
   \qed

 \section{Concluding Remark}
Our results imply a tight connection between $C_k$-rainbow colorings in the cube and constructions of large generalized Sidon sets. When $k \equiv 0$ (mod 4) Construction 1 gives the correct order of magnitude, however for $k \equiv 2$ (mod 4) the same method does not work. In this case an approach similar to Construction 2 would work provided we can construct large sets that do not contains solutions to certain equations.

\begin{conj}  \label{conj}  Fix $4 \le k \equiv 2$ {\rm (mod 4)}.  Then $f(n,k)=n^{\lfloor k/4 \rfloor + o(1)}$.
\end{conj}
For the first open case $k=10$, we can show that $f(n,10)=n^{2+o(1)}$ provided one can construct a set $S \subset [N]$ with $|S| > N^{1/2-o(1)}$ that contains no nontrivial solution to any of the following equations:
$$x_1+ x_2= x_3+x_4$$
$$x_1+x_2+x_3= x_4+2x_5$$
$$x_1 + 2x_2= x_3+2x_4.$$

Ruzsa~\cite{R1, R2}  defined the genus $g(E)$ of an equation 
$$ E: \quad a_1x_1 + \cdots  + a_kx_k = 0$$
 as the largest $m$ such that there is a partition  $S_1\cup \ldots \cup S_m$ of $[k]$ where the $S_i$ are disjoint, non-empty  and for all $j$, 
\begin{equation} \label{re} \sum\limits_{i \in S_j} a_i = 0. \end{equation}
 A solution $(x_1, \ldots, x_k)$ of $E$ is  trivial if  there are $l$ distinct numbers among $\{x_1, \ldots, x_k\}$ and  (\ref{re}) holds for a partition $S_1 \cup \ldots  \cup S_l$ of $[k]$ into disjoint, non-empty parts such that $x_i = x_j$ if and only if $i,j \in S_v$ for some $v$.

Ruzsa showed that if  $S \subset [n]$ has no nontrivial solutions to $E$ then $|S|\le O(n^{1/g(E)})$. The question of whether there exists $S$ with 
$|S|=n^{1/g(E)-o(1)}$ remains open for most equations $E$. The set of equations above has genus two so it is plausible
  that one can prove Conjecture \ref{conj} for $k=10$ using this approach.

For the general case, we can provide a rainbow coloring if our set $S$ contains no nontrivial solutions to any of the three equations below  with $m = \lfloor k/4 \rfloor$.

\begin{center}
$x_1 + \cdots  + x_m  = x_{m+1} + \cdots  + x_{2m}$
\\
$x_1 + \cdots  + x_m + x_{m+1} = x_{m+2} + \cdots  + x_{2m} + 2x_{2m+1}$
\\
$x_1 + \cdots  + x_{m-1} + 2x_m = x_{m+1} + \cdots  + x_{2m-1} + 2x_{2m}$.
\end{center}

The set  of equations above 
has genus $m = \lfloor k/4 \rfloor$, so if Ruzsa's question has a positive 
answer, then we would be able to construct a set of the desired size.

\end{document}